\documentclass[10pt]{amsart}
\usepackage{mathrsfs}
\usepackage{amsfonts} 
\usepackage{graphicx,latexsym,bm,amsmath,amssymb,verbatim,multicol,lscape}
\theoremstyle{definition}

\theoremstyle{remark}

\numberwithin{equation}{section}

\begin{document}
\title[LCM of consecutive arithmetic progression terms]
{The least common multiple of consecutive arithmetic progression terms}%
\author{Shaofang Hong}
\address{Mathematical College, Sichuan University, Chengdu 610064, P.R. China}
\email{sfhong@scu.edu.cn, sfhong09@gmail.com, s-f.hong@tom.com,
hongsf02@yahoo.com }
\author{Guoyou Qian}
\address{Mathematical College, Sichuan University, Chengdu 610064, P.R. China}
\email{qiangy1230@gmail.com, qiangy1230@163.com, qiangy0901@tom.com}
\thanks{The research was supported partially by National Science Foundation of China
Grant \# 10971145 and by Program for New Century Excellent Talents
in University Grant \# NCET-06-0785}

\keywords{arithmetic progression, least common multiple, $p$-adic
valuation, arithmetic function, the smallest period,}
\subjclass[2000]{Primary 11B25, 11N13, 11A05}
\date{\today}%
\begin{abstract}
Let $k\ge 0,a\ge 1$ and $b\ge 0$ be integers. We define the
arithmetic function $g_{k,a,b}$ for any positive integer $n$ by
$g_{k,a,b}(n):=\frac{(b+na)(b+(n+1)a)\cdots (b+(n+k)a)} {{\rm
lcm}(b+na,b+(n+1)a,\cdots,b+(n+k)a)}.$ Letting $a=1$ and $b=0$, then
$g_{k,a,b}$ becomes the arithmetic function introduced previously by
Farhi. Farhi proved that $g_{k,1,0}$ is periodic and that $k!$ is a
period. Hong and Yang improved Farhi's period $k!$ to ${\rm
lcm}(1,2,...,k)$ and conjectured that $\frac{{\rm
lcm}(1,2,...,k,k+1)}{k+1}$ divides the smallest period of
$g_{k,1,0}$. Recently, Farhi and Kane proved this conjecture and
determined the smallest period of $g_{k,1,0}$. For the general
integers $a\ge 1$ and $b\ge 0$, it is natural to ask the interesting
question: Is $g_{k,a,b}$ periodic? If so, then what is the smallest
period of $g_{k,a,b}$? We first show that the arithmetic function
$g_{k,a,b}$ is periodic. Subsequently, we provide detailed $p$-adic
analysis of the periodic function $g_{k,a,b}$. Finally, we determine
the smallest period of $g_{k,a,b}$. Our result extends the
Farhi-Kane theorem from the set of positive integers to general
arithmetic progressions.
\end{abstract}

\maketitle

\section{\bf Introduction}

Many beautiful and important theorems about the arithmetic
progression in number theory are known. Dirichlet's theorem
\cite{[A]} \cite{[HW]} and the Green-Tao theorem \cite{[GT]} are the
two most famous examples. For some other results, see, for example,
 \cite{[BBGH]}, \cite{[HLe]}, \cite{[HLo]}, \cite{[S]} and \cite{[SS]}.
Meanwhile, the least common multiple of any given sequence of
positive integers has received lots of attentions by many authors,
see, for instance, \cite{[BK]}, \cite{[BKS]},
\cite{[C]}-\cite{[F2]}, \cite{[Ha]}-\cite{[HK]}, \cite{[HY1]},
\cite{[MS]} and \cite{[N]}. For the detailed background information
about the least common multiple of finite arithmetic progressions,
we refer the readers to \cite{[HY2]}.

Farhi \cite{[F1]}, \cite{[F2]} investigated the least common
multiple of a finite number of consecutive integers. Let $k\ge 0$ be
an integer. It was proved in \cite{[F1]} and \cite{[F2]} that
lcm$(n, n+1, ..., n+k)$ is divisible by $n{n+k \choose k}$ and also
divides $n{n+k \choose k}$lcm$({k \choose 0}, {k \choose 1}, ..., {k
\choose k})$. Farhi \cite{[F1]}, \cite{[F2]} showed that the last
equality holds if $k!|(n+k+1)$. Farhi also introduced the arithmetic
function $g_k$, which is defined for any positive integer $n$ by
$$
g_k(n):=\frac{n(n+1)\cdots(n+k)}{{\rm lcm}(n,n+1,\cdots,n+k)}.
$$
Farhi then proved that the sequence $\{g_k\}_{k=0}^\infty $
satisfies the recursive relation $g_k(n)={\rm
gcd}(k!,(n+k)g_{k-1}(n))$ for all positive integers $n$, where
gcd$(a, b)$ means the greatest common divisor of integers $a$ and
$b$. Using this relation, we can easily show (by induction on $k$)
that for any nonnegative integer $k$, the function $g_k$ is periodic
of period $k!$. This is a result due to Farhi \cite{[F2]}. Define
$P_k$ to be the smallest period of the function $g_k$. Then Farhi's
result says that $P_k|k!$. Define $L_0:=1$ and for any integer $k\ge
1$, we define $L_k:={\rm lcm}(1,2,...,k)$. Hong and Yang
\cite{[HY2]} showed that $g_k(1)|g_k(n)$ for any nonnegative integer
$k$ and any positive integer $n$. Consequently, using this result,
they showed that $P_k|L_k$ for all positive integer $k$. This
improves Farhi's period. In \cite{[HY2]}, Hong and Yang raised a
conjecture stating that $\frac{L_{k+1}}{k+1}$ divides $P_k$ for all
nonnegative integers $k$. From this conjecture, one can read that
$k|P_k$ and $P_k=L_k$ if $k+1$ is a prime. Very recently, Farhi and
Kane \cite{[FK]} found a proof of the Hong-Yang conjecture.
Furthermore, Farhi and Kane determined the exact value of $P_k$,
which solved the open problem posed by Farhi in \cite{[F2]}.

Throughout this paper, let ${\mathbb Q}$ and ${\mathbb N}$ denote
the field of rational numbers and the set of positive integers.
Define ${\mathbb N}_0:={\mathbb N}\cup\{0\}$. Let $k,b\in
\mathbb{N}_0$ and $a\in \mathbb{N}$. We define the arithmetic
function $g_{k,a,b}:\mathbb{N}\longrightarrow \mathbb{N}$ by
$$
g_{k,a,b}(n)=\frac{(b+na)(b+(n+1)a)\cdots (b+(n+k)a)} {{\rm
lcm}(b+na,b+(n+1)a,\cdots,b+(n+k)a)}.
$$
Note that $g_{k,1,0}=g_k$. It is
natural to ask the following interesting question:\\

\noindent{\bf Problem 1.1.} Let $k\ge 0,a\ge 1$ and $b\ge 0$ be
integers. Is $g_{k,a,b}$ periodic and, if so, what is the smallest
period of $g_{k,a,b}$?\\

Assume that $g_{k,a,b}$ is periodic and $P_{k,a,b}$ is the smallest
period of $g_{k,a,b}$. Then we can use $P_{k,a,b}$ to give a formula
for ${\rm lcm}(b+na,b+(n+1)a,\cdots,b+(n+k)a)$ as follows: For any
positive integer $n$, we have
$$
{\rm lcm}(b+na,b+(n+1)a,\cdots,b+(n+k)a)
=\frac{(b+na)(b+(n+1)a)\cdots (b+(n+k)a)} {g_{k,a,b}(\langle
n\rangle_{P_{k,a,b}})},
$$
where $\langle n\rangle_{P_{k,a,b}}$ means the least nonnegative
residue of $n$ modulo $P_{k,a,b}$. Therefore it is important to
determine the exact value of $P_{k,a,b}$.

In this paper, we investigate the least common multiple of
consecutive terms in arithmetic progressions. As usual, for any
prime number $p$, we let $v_{p}$ be the normalized $p$-adic
valuation of ${\mathbb Q}$, i.e., $v_p(a)=s$ if $p^{s}\parallel a$.
For any real number $x$, by $\lfloor x\rfloor$ we denote the largest
integer no more than $x$. Let $e_{p,k}:=\lfloor {\rm
log}_{p}k\rfloor ={\rm max}_{1\leq i\leq k}\{v_{p}(i)\}$ be the
largest exponent of a power of $p$ that is at
most $k$. We can now give the main result of this paper as follows.\\

\noindent{\bf Theorem 1.2.} {\it Let $k\ge 0, a\ge 1$ and $b\ge 0$
be integers. Then the arithmetic function $g_{k,a,b}$ is periodic,
and if ${\rm gcd}(a, b)=1$, then its smallest period equals $Q_{k,
a}$, where
$$
Q_{k,a}:=\frac{L_k}{\delta _{k,a}\cdot \displaystyle\prod_{{\rm
prime} \ q|{\rm gcd}(a, L_k)}q^{e_{q,k}}}, \eqno(1.1)
$$
and
\begin{align*}
\delta _{k,a}:={\left\{
  \begin{array}{rl}
p^{e_{p,k}}, &\text{if} \ p\nmid a \
\text{and} \ v_p(k+1)\geq e_{p,k} \ \text{for some prime} \ p\le k,\\
1, &\text{otherwise.}
 \end{array}
\right.}
\end{align*}
For ${\rm gcd}(a, b)>1$, its smallest period is equal to $Q_{k, a'}$
with
$a'=\frac{a}{{\rm gcd}(a, b)}$.}\\
\\
Thus we answer completely Problem 1.1. Our result extends the
Farhi-Kane theorem from the set of positive integers to general
arithmetic progressions.

The paper is organized as follows. In Section 2, by using a
well-known result of Hua \cite{[Hu]}, we show that the arithmetic
function $g_{k,a,b}$ is periodic; see Theorem 2.5. Then in Section
3, we provide detailed $p$-adic analysis to the periodic function
$g_{k,a,b}$ and determine the smallest period of $g_{k,a,b}$. In the
last section, we show Theorem 1.2 and give an example to illustrate
the validity of Theorem 1.2.

\section{\bf The periodicity of $g_{k,a,b}$}

Hong and Yang \cite{[HY2]} proved that $L_k$ is a period of $g_{k}$.
In this section, we introduce a new method to show that for any
integers $k\ge 0, a\ge 1$ and $b\ge 0$, the arithmetic function
$g_{k,a,b}$ is periodic, and particularly $L_k$ is also a period of
$g_{k,a,b}$. First we need a well-known result of Hua. One can
easily deduce this result by using the principle of
inclusion-exclusion (see, for instance, Page 11 of \cite{[Hu]}).\\

\noindent{\bf Lemma 2.1.} \cite{[Hu]} {\it Let $a_1,a_2,\ldots,a_n$
be any $n$ positive integers. Then we have}
$$
{\rm lcm}(a_1,a_2,\ldots,a_n)=a_1a_2\ldots a_n\cdot
\prod_{r=2}^{n}\prod_{1\leq i_1<...<i_{r}\leq n}({\rm gcd}(a_{i_1},
..., a_{i_{r}}))^{(-1)^{r-1}}.
$$\\

\noindent{\bf Lemma 2.2.} {\it Let $a_1,a_2,\ldots,a_n$ and
$b_1,b_2,\ldots,b_n$ be any $2n$ positive integers. Let $3\le t\le
n$ be a given integer. If ${\rm gcd}(a_{i_1}, ..., a_{i_t})={\rm
gcd}(b_{i_1}, ..., b_{i_t})$ for any $1\leq i_1<...<i_t\leq n$, then
we have}
$$
\frac{a_1a_2\cdots a_n}{{\rm lcm}(a_1,a_2,\ldots,a_n)}\cdot
\prod_{r=2}^{t-1}\prod_{1\leq i_1<...<i_{r}\leq n}({\rm
gcd}(a_{i_1}, ..., a_{i_{r}}))^{(-1)^{r-1}}$$
$$=\frac{b_1b_2\cdots b_n}{{\rm
lcm}(b_1,b_2,\ldots,b_n)}\cdot \prod_{r=2}^{t-1}\prod_{1\leq
i_1<...<i_{r}\leq n}({\rm gcd}(b_{i_1}, ...,
b_{i_{r}}))^{(-1)^{r-1}}.
$$

\begin{proof}
If ${\rm gcd}(a_{i_1}, ..., a_{i_t})={\rm gcd}(b_{i_1}, ...,
b_{i_t})$ for any $1\leq i_1<...<i_t\leq n$, then we have ${\rm
gcd}(a_{i_1}, ..., a_{i_k})={\rm gcd}(b_{i_1}, ..., b_{i_k})$ for
any $1\leq i_1<...<i_k\leq n$ and any $n\geq k\geq t$. Thus, by
using Lemma 2.1, we get the result of Lemma 2.2.
\end{proof}

In particular, we have the following result.\\

\noindent{\bf Lemma 2.3.} {\it Let $a_1,a_2,\ldots,a_n$ and
$b_1,b_2,\ldots,b_n$ be any $2n$ positive integers. If for any
$1\leq i_1<i_2<i_3\leq n$, we have ${\rm gcd}(a_{i_1}, a_{i_2},
a_{i_3})={\rm gcd}(b_{i_1}, b_{i_2}, b_{i_3})$, then
\begin{align*}
\frac{1}{\prod_{1\leq i<j\leq n}{\rm gcd}(a_{i},a_{j})}\cdot
\frac{a_1a_2\cdots a_n}{{\rm lcm}(a_1,a_2,\ldots,a_n)}
=\frac{1}{\prod_{1\leq i<j\leq n}{\rm gcd}(b_{i},b_{j})}\cdot
\frac{b_1b_2\cdots b_n}{{\rm lcm}(b_1,b_2,\ldots,b_n)}.
\end{align*}
}
\begin{proof}
Since ${\rm gcd}(a_{i_1},a_{i_2},a_{i_3})={\rm
gcd}(b_{i_1},b_{i_2},b_{i_3})$ for any $1\leq i_1<i_2<i_3\leq n$, so
we have ${\rm gcd}(a_{i_1},\ldots,a_{i_k})={\rm
gcd}(a_{i_1},\ldots,a_{i_k})$ for any $1\leq i_1<\cdots<i_k\leq n$
and $ k\geq 3$. By using Lemma 2.1, we get the conclusion of Lemma
2.3.
\end{proof}

Notice that if ${\rm gcd}(a_{i}, a_{j})={\rm gcd}(b_{i }, b_{j})$
for any $1\leq i<j\le n$, then ${\rm
gcd}(a_{i_1},a_{i_2},a_{i_3})={\rm gcd}(b_{i_1},b_{i_2},b_{i_3})$
for any $1\leq i_1<i_2<i_3\leq n$.
It follows immediately from Lemma 2.3 that the following is true.\\

\noindent{\bf Corollary 2.4.} {\it Let $a_1,a_2,\ldots,a_n$ and
$b_1,b_2,\ldots,b_n$ be any $2n$ positive integers. If ${\rm
gcd}(a_{i}, a_{j})={\rm gcd}(b_{i }, b_{j})$ for any $1\leq i<j\leq
n$, then we have
\begin{align*}
\frac{a_1a_2\cdots a_n}{{\rm lcm}(a_1,a_2,\ldots,a_n)}
=\frac{b_1b_2\cdots b_n}{{\rm lcm}(b_1,b_2,\ldots,b_n)}.
\end{align*}}

We can now give the main result of this section. This also gives an
alternative proof to the Hong-Yang period of the
periodic function $g_k$ \cite{[HY2]}.\\

\noindent{\bf Theorem 2.5.} {\it Let $k\ge 0,a\ge 1$ and $b\ge 0$ be
integers. Then the arithmetic function $g_{k,a,b}$ is periodic, and
$L_k$ is a period of $g_{k,a,b}$. }

\begin{proof}
Let $n$ be any positive integer. For any $0\leq i<j\leq k$, we have
\begin{align*}
{\rm gcd}(b+(n+i+L_k)a,b+(n+j+L_k)a)&={\rm gcd}(b+(n+i+L_k)a,(j-i)a)\\
&={\rm gcd}(b+(n+i)a,(j-i)a)\\
&={\rm gcd}(b+(n+i)a,b+(n+j)a).
\end{align*}
Thus by Corollary 2.4, we obtain
\begin{align*}
&\frac{(b+(n+L_k)a)(b+(n+1+L_k)a)\cdots (b+(n+k+L_k)a)}{{\rm
lcm}(b+(n+L_k)a,b+(n+1+L_k)a,\cdots,b+(n+k+L_k)a) }\\
&=\frac{(b+na)(b+(n+1)a)\cdots (b+(n+k)a)}{{\rm
lcm}(b+na,b+(n+1)a,\cdots,b+(n+k)a)}.
\end{align*}
In other words, for any positive integer $n$, we have
$g_{k,a,b}(n+L_k)=g_{k,a,b}(n)$ as desired.
\end{proof}

Evidently, Theorem 2.5 gives an affirmative answer to the first part
of Problem 1.1.

\section{\bf $p$-Adic analysis of $g_{k,a,b}$}

Throughout this section, we always let $k\ge 0,a\ge 1$ and $b\ge 0$
be integers such that gcd$(a, b)=1$. From the main result of
previous section (Theorem 2.5), we know that the arithmetic function
$g_{k,a,b}$ is periodic. Let $P_{k,a,b}$ denote the smallest period
of $g_{k,a,b}$. Then by Theorem 2.5 we know that $P_{k,a,b}$ is a
divisor of $L_k$. But the exact value of $P_{k,a,b}$ is still
unknown. In this section, we will determine the exact value of
$P_{k,a,b}$. We need some more notation. Let
$$S_{k,a,b}(n):=\{b+na, b+(n+1)a, ..., b+(n+k)a\}$$
be any $k+1$ consecutive terms in the arithmetic progression
$\{b+ma\}_{m\in \mathbb{N}_0}$. For a given prime number $p$, define
$g_{p,k,a,b}(n):=v_{p}(g_{k,a,b}(n))$. Since $g_{k,a,b}$ is a
periodic function, $g_{p,k,a,b}$ is also a periodic function for
each prime $p$ and $P_{k,a,b}$ is a period of $g_{p,k,a,b}$. Let
$P_{p,k,a,b}$ be the smallest period of $g_{p,k,a,b}$.
We have the following result.\\

\noindent{\bf Lemma 3.1.} {\it We have $ P_{k,a,b}=\displaystyle{\rm
lcm}_{p \ {\rm prime}}(P_{p,k,a,b}).$}

\begin{proof}
Since for any $n\in \mathbb{N}$, we have that
$v_{p}(g_{k,a,b}(n+P_{k,a,b}))=v_{p}(g_{k,a,b}(n))$, i.e.,
$P_{p,k,a,b}\mid P_{k,a,b}$ for each prime $p$. Hence we have ${\rm
lcm}_{p \ {\rm prime}}(P_{p,k,a,b})\mid P_{k,a,b}$.

Conversely, for any $n\in \mathbb{N}$, we have that
$v_{p}(g_{k,a,b}(n+{\rm lcm}_{p \ {\rm
prime}}(P_{p,k,a,b}))=v_{p}(g_{k,a,b}(n))$ for each prime $p$. Thus,
we have $g_{k,a,b}(n+{\rm lcm}_{p \ {\rm
prime}}(P_{p,k,a,b}))=g_{k,a,b}(n)$ for any $n\in \mathbb{N}$, that
is, we have $P_{k,a,b}\mid {\rm lcm}_{p \ {\rm
prime}}(P_{p,k,a,b})$. Therefore, we have $P_{k,a,b}={\rm lcm}_{p \
{\rm prime}}(P_{p,k,a,b})$, as required.
\end{proof}

Hence we only need to compute $P_{p,k,a,b}$ for each prime $p$ to
get the exact value of $P_{k,a,b}$. The following result is due to
Farhi \cite{[F1]}. An alternative proof of it was given by Hong
and Feng \cite{[HF]}.\\

\noindent{\bf Lemma 3.2.} \cite{[F1]} \cite{[HF]} {\it Let
$\{u_i\}_{i\in \mathbb{N}_0}$ be a strictly increasing arithmetic
progression of non-zero integers and $k$ be any given non-negative
integer. Then the integer ${\rm lcm}(u_0,u_1, ...,\\u_k)$ is a
multiple of
$\frac{u_0u_1\cdots u_k}{k!({\rm gcd}(u_0,u_1))^k}$.}\\

\noindent{\bf Lemma 3.3.} {\it For any positive integer $n$, we have
$g_{k,a,b}(n)\mid k!.$}

\begin{proof}
Let $u_i=b+a(n+i)$ for $0\le i\le k$. Then ${\rm gcd}(u_0,u_1)=1$
since $a$ and $b$ are coprime. So by Lemma 3.2 we know that there is
an integer $A$ such that
$$
{\rm lcm}(b+na, b+(n+1)a, ..., b+(n+k)a)=A\cdot
\frac{(b+an)(b+a(n+1))\cdots (b+a(n+k)}{k!}.
$$
It then follows that $k!=A\cdot g_{k,a,b}(n)$.
\end{proof}

It follows from Lemma 3.3 that $g_{p,k,a,b}(n)=v_{p}(g_{k, a,
b}(n))= 0$ for each prime $p>k$ and any positive integer $n$. Hence
we have $P_{p, k, a, b}=1$ for each prime $p>k$. So by Lemma 3.1, in
order to determine the exact value of $P_{k,a,b}$, it suffices to
compute the exact value of $P_{p,k,a,b}$ for all the primes $p$ such
that $1<p\leq k$. First we consider the case that $p\mid a$ and
$1<p\leq k$. Since ${\rm gcd}(a,b)=1$, we have gcd$(p,b)=1$, and
thus ${\rm gcd}(p, b+(n+i)a)=1$ for any integers $n\in \mathbb{N}$
and $0\leq i\leq k$. Hence ${\rm gcd}(p,g_{k,a,b}(n))=1$ for any
integer $n\ge 1$, i.e., we have $g_{p,k,a,b}(n)=0$ for any integer
$n\ge 1$ if $p\mid a$. Thus $P_{p,k,a,b}=1$ if $p\mid a$. We put
these facts
into the following lemma.\\

\noindent{\bf Lemma 3.4.} {\it Let $p$ be a prime such that either
$p>k$ or $p|a$. Then we have $P_{p,k,a,b}=1$.}\\

In what follows we treat the remaining case that $p\nmid a$ and
$1<p\leq k$. Clearly we have

\begin{align*}
g_{p,k,a,b}(n)&= \sum_{m\in S_{k,a,b}(n)}v_{p}(m)-
{\rm max}_{m\in S_{k,a,b}(n)}v_{p}(m)\\
&= \sum_{e\ge 1}\sum_{m\in S_{k,a,b}(n)}(1 \ {\rm if} \ p^e\mid
m)-\sum_{e\ge 1}(1 \ {\rm if} \ p^e \
{\rm divides \ some} \ m\in S_{k,a,b}(n))\\
&= \sum_{e\ge 1}\#\{m\in S_{k,a,b}(n):p^e\mid m\}-\sum_{e\ge 1}(1 \
{\rm if} \ p^e \
{\rm divides \ some} \ m\in S_{k,a,b}(n))\\
&= \sum_{e\ge 1}{\rm max} (0, \#\{m\in S_{k,a,b}(n):p^e\mid m\}-1).
\ \ \ \ \ \ \ \ \ \ \ \ \ \ \ \ \ \ \ \ \ \ \ \ \ \ \ \ \ \ \ \ \
(3.1)
\end{align*}
Then we have the following lemmas:\\

\noindent{\bf Lemma 3.5.} {\it If $p\nmid a$ and $e>e_{p,k}$, then
there is at most one element of $S_{k,a,b}(n)$ which is divisible by
$p^{e}$.}

\begin{proof}
Suppose that there exist two integers $i,j$ such that $p^{e}\mid
b+(n+i)a$ and $p^{e}\mid b+(n+j)a$, where $0\leq i<j\leq k$, then we
have $p^{e}\mid (j-i)a$. Since ${\rm gcd}(p,a)=1$, we get $p^{e}\mid
(j-i)$. From it we deduce that $v_p(j-i)\ge e>e_{p,k}$. This is a
contradiction.
\end{proof}

\noindent{\bf Lemma 3.6.} {\it Let $e$ be a positive integer. If
$p\nmid a$, then any $p^e$ consecutive terms in the arithmetic
progression $\{b+ma\}_{m\in \mathbb{N}_0}$ are pairwise incongruent
modulo $p^e$. Furthermore, if $e\leq e_{p,k}$, then there is at
least one element of $S_{k,a,b}(n)$ divisible by $p^{e}$ .}

\begin{proof}
Suppose that there exist two integers $i,j$ such that
$b+(m+i)a\equiv b+(m+j)a \pmod {p^e}$, where $m\ge 0$ and $0\leq
i<j\leq p^e-1$. Then $p^e\mid (j-i)a$. Since ${\rm gcd}(p,a)=1$, we
have $p^{e}\mid (j-i)$. This is impossible. Thus the first part is
true.

Now let $e\leq e_{p,k}$. Then $1\leq p^{e}\leq k$. Hence
$S_{k,a,b}(n)$ holds $p^e$ consecutive terms and one of which is
divisible by $p^e$ by the above discussion. So the second part
holds.
\end{proof}

By Lemma 3.5, we know that all the terms in the right-hand side of
(3.1) are 0 if $e>e_{p,k}$. By Lemma 3.6, there is at least one
element divisible by $p^e$ in the set $S_{k,a,b}(n)$ if $e\leq
e_{p,k}$. Therefore we obtain by (3.1)
$$g_{p,k,a,b}(n)=\sum_{e=1}^{e_{p,k}}f_{e}(n), \eqno (3.2)
$$
where $f_{e}(n):=\#\{m\in S_{k,a,b}(n):p^e\mid m\}-1$. Since
$b+(n+i+p^e)a\equiv b+(n+i)a\pmod {p^e}$ for any $i\in
\{0,1,\ldots,k\}$, we have $f_{e}(n+p^e)=f_{e}(n)$. Therefore $p^e$
is a period of $f_{e}(n)$. Hence $f_{e}(n+p^{e_{p,k}})=f_{e}(n)$ is
true for each $e\in \{1,\ldots,e_{p,k}\}$. This implies that
$g_{p,k,a,b}(n+p^{e_{p,k}})=g_{p,k,a,b}(n)$. Consequently,
$p^{e_{p,k}}$ is a period of $g_{p,k,a,b}(n)$. Thus $P_{p,k,a,b}\mid
p^{e_{p,k}}$. It follows immediately that the $P_{p,k,a,b}$ are
relatively prime for different prime numbers $p$. But Lemma 3.1 and
Lemma 3.4 tell us that $P_{k,a,b}=\displaystyle{\rm lcm}_{p \ {\rm
prime}, p\le k, p\nmid a}(P_{p,k,a,b})$. Therefore we get the following result.\\

\noindent{\bf Lemma 3.7.} {\it We have
$$P_{k,a,b}=\prod_{p \ {\rm prime}, p\nmid a, p\leq k}P_{p,k,a,b},
$$
where $P_{p,k,a,b}$ satisfies that $P_{p,k,a,b}|p^{e_{p,k}}$.}\\

According to Lemma 3.7, it suffices to compute the $p$-adic
valuation of $P_{p,k,a,b}$ for the prime numbers $p$ satisfying
$p\nmid a$ and $p\in (1,k]$. Now let us determine the $p$-adic
valuation of $P_{k,a,b}$ for these prime numbers $p$.\\

\noindent{\bf Proposition 3.8.} {\it Let $a\ge 1$ and $b\ge 0$ be
integers such that ${\rm gcd}(a,b)=1$. Let $k\geq 2$ be an integer
and $p\in (1,k]$ be a prime number such that $p\nmid a$.

{\rm (i).} If $v_{p}(k+1)<e_{p,k}$, then $v_{p}(P_{k,a,b})=e_{p,k}.$

{\rm (ii).} If $v_{p}(k+1)\geq e_{p,k},$ then $v_{p}(P_{k,a,b})=0.$}

\begin{proof}
(i). Since $p^{e_{p,k}}$ is a period of $g_{p,k,a,b}$, it suffices
to prove that $p^{e_{p,k}-1}$ is not the period of $g_{p,k,a,b}$,
from which it follows that $p^{e_{p,k}}$ is the smallest period of
$g_{p,k,a,b}$. By (3.2), we have
$$
g_{p,k,a,b}(n)=\sum_{e=1}^{e_{p,k}}f_{e}(n)=\sum_{e=1}^{e_{p,k}-1}f_{e}(n)+f_{e_{p,k}}(n).
$$
Since $p^{e_{p,k}-1}$ is a period of
$\sum_{e=1}^{e_{p,k}-1}f_{e}(n)$, it is sufficient to prove that
$p^{e_{p,k}-1}$ is not the period of $f_{e_{p,k}}(n)$. We claim that
there exists a positive integer $n_0$ such that
$f_{e_{p,k}}(n_0+p^{e_{p,k}-1})=f_{e_{p,k}}(n_0)-1$.

By $v_{p}(k+1)<e_{p,k}$, we deduce that $p^{e_{p,k}}\nmid (k+1)$ and
$p^{e_{p,k}}\le k$. We can suppose that $k+1\equiv l\pmod
{p^{e_{p,k}}} \ {\rm for \ some } \ 1\leq l\leq p^{e_{p,k}}-1.$ We
divide the proof of part (i) into the following two cases:

{\sc Case 1.} $1\leq l\leq p^{e_{p,k}}-p^{e_{p,k}-1}$. Since $p\nmid
a$, we can always find a suitable $n_0$ such that $b+n_{0}a\equiv
0\pmod {p^{e_{p,k}}}$. Consider the following two sets:
$$S_{k,a,b}(n_0)=\{b+n_{0}a,\ldots,b+(n_{0}+p^{e_{p,k}-1}-1)a,
b+(n_{0}+p^{e_{p,k}-1})a,\ldots,b+(n_{0}+k)a\}$$
and
\begin{align*}
S_{k,a,b}(n_0+p^{e_{p,k}-1})=\{& b+(n_0+p^{e_{p,k}-1})a,\ldots,
b+(n_0+k)a,\\
&b+(n_0+k+1)a,\ldots, b+(n_0+k+p^{e_{p,k}-1})a\}.
\end{align*}
Now $\{b+(n_0+p^{e_{p,k}-1})a,\ldots, b+(n_0+k)a\}$ is the
intersection of $S_{k,a,b}(n_0)$ and $S_{k,a,b}(n_0+p^{e_{p,k}-1})$.
So to compare the number of terms divisible by $p^{e_{p,k}}$ in the
set $S_{k,a,b}(n_0)$ with the number of terms divisible by
$p^{e_{p,k}}$ in the set $S_{k,a,b}(n_0+p^{e_{p,k}-1})$, it suffices
to compare the number of terms divisible by $p^{e_{p,k}}$ in the set
$\{b+n_{0}a,\ldots,b+(n_0+p^{e_{p,k}-1}-1)a\}$ with the number of
terms divisible by $p^{e_{p,k}}$ in the set
$\{b+(n_{0}+k+1)a,\ldots,b+(n_0+k+p^{e_{p,k}-1})a\}$. By Lemma 3.6,
any $p^{e_{p,k}}$ consecutive terms in the arithmetic progression
$\{b+ma\}_{m\in \mathbb{N}_0}$ are pairwise incongruent modulo
$p^{e_{p,k}}$. Thus the terms divisible by $p^{e_{p,k}}$ in the
arithmetic progression $\{b+ma\}_{m\in \mathbb{N}_0}$ must be of the
form $b+(n_0+tp^{e_{p,k}})a, \ t\in \mathbb{Z}$. Since $k+1\equiv
l\pmod {p^{e_{p,k}}}$ and $1\leq l\leq p^{e_{p,k}}-p^{e_{p,k}-1}$,
we have $k+j\equiv l+j-1\not\equiv 0 \pmod {p^{e_{p,k}}}$ for all
$1\leq j\leq p^{e_{p,k}-1})$. Hence $p^{e_{p,k}}\nmid (b+(n_0+k+j)a$
for all $1\leq j\leq p^{e_{p,k}-1}$. Thus all the elements in the
set $\{b+(n_0+k+1)a,\ldots,b+(n_0+k+p^{e_{p,k}-1})a\}$ are not
divisible by $p^{e_{p,k}}$. On the other hand, since $b+an_0\equiv
0\pmod {p^{e_{p,k}}}$, it follows from Lemma 3.6 that there is
exactly one term in the set
$\{b+n_{0}a,b+(n_{0}+1)a,\ldots,b+(n_{0}+p^{e_{p,k}-1}-1)a\}$ which
is divisible by $p^{e_{p,k}}$. Therefore the number of terms
divisible by $p^{e_{p,k}}$ in the set
$S_{k,a,b}(n_{0}+p^{e_{p,k}-1})$ is equal to the number of terms
divisible by $p^{e_{p,k}}$ in the set $S_{k,a,b}(n_0)$ minus one.
Namely, $f_{e_{p,k}}(n_0+p^{e_{p,k}-1})=f_{e_{p,k}}(n_0)-1$ as
required. The claim is proved in this case.

{\sc Case 2}. $p^{e_{p,k}}-p^{e_{p,k}-1}<l\leq p^{e_{p,k}}-1$. Since
$p\nmid a$, it is easy to see that there is a positive integer $n_0$
such that $b+(n_0+p^{e_{p,k}-1}-1)a\equiv 0\pmod {p^{e_{p,k}}}$. As
in the discussion of Case 1, to compare the number of terms
divisible by $p^{e_{p,k}}$ in the set $S_{k,a,b}(n_0)$ with the
number of terms divisible by $p^{e_{p,k}}$ in the set
$S_{k,a,b}(n_0+p^{e_{p,k}-1})$, it suffices to compare the number of
terms divisible by $p^{e_{p,k}}$ in the set
$\{b+n_{0}a,\ldots,b+(n_0+p^{e_{p,k}-1}-1)a\}$ with the number of
terms divisible by $p^{e_{p,k}}$ in the set
$\{b+(n_{0}+k+1)a,\ldots,b+(n_0+k+p^{e_{p,k}-1})a\}$. From
$b+(n_0+p^{e_{p,k}-1}-1)a\equiv 0\pmod {p^{e_{p,k}}}$ one can deduce
that the terms divisible by $p^{e_{p,k}}$ in the arithmetic
progression $\{b+ma\}_{m\in \mathbb{N}_0}$ must be of the form
$b+(n_0+p^{e_{p,k}-1}-1+tp^{e_{p,k}})a$ with $t\in \mathbb{Z}$.
Since $k+1\equiv l\pmod {p^{e_{p,k}}}$ for some
$p^{e_{p,k}}-p^{e_{p,k}-1}<l\leq p^{e_{p,k}}-1$, we have
$p^{e_{p,k}}-p^{e_{p,k}-1}+1\leq l+j-1\le
p^{e_{p,k}}+p^{e_{p,k}-1}-2$ and so $k+j\equiv l+j-1\not\equiv
p^{e_{p,k}-1}-1\pmod {p^{e_{p,k}}}$ for all $1\leq j\leq
p^{e_{p,k}-1}$. It follows that for all $1\leq j\leq p^{e_{p,k}-1}$,
we have $p^{e_{p,k}}\nmid (b+(n_0+k+j)a)$. That is, there does not
exist an integer divisible by $p^{e_{p,k}}$ in the set
$\{b+(n_0+k+1)a,\ldots,b+(n_0+k+p^{e_{p,k}-1})a\}$. But the term
$b+(n_0+p^{e_{p,k}-1}-1)a$ is the only term divisible by
$p^{e_{p,k}}$ in the set
$\{b+n_{0}a,b+(n_0+1)a,\ldots,b+(n_0+p^{e_{p,k}-1}-1)a\}$. Thus the
number of terms divisible by $p^{e_{p,k}}$ in the set
$S_{k,a,b}(n_0+p^{e_{p,k}-1})$ equals the number of terms divisible
by $p^{e_{p,k}}$ in the set $S_{k,a,b}(n_0)$ minus one. Hence the
desired result $f_{e_{p,k}}(n_0+p^{e_{p,k}-1})=f_{e_{p,k}}(n_0)-1$
follows immediately. The proof of the claim is complete.

From the claim we deduce immediately that $p^{e_{p,k}-1}$ is not a
period of $g_{p,k,a,b}$. Thus $p^{e_{p,k}}$ is the smallest period
of $g_{p,k,a,b}$. It follows that $v_{p}(P_{k,a,b})=e_{p,k}$ as
desired.

(ii). By Lemma 3.7, we know that to prove part (ii), it is
sufficient to prove that $v_p(P_{q,k,a,b})=0$ for each prime $q$
with $q\le k$ and $q\nmid a$. For any prime $q$ different from $p$,
since $P_{q,k,a,b}|q^{e_{q,k}}$, we then have $v_p(P_{q,k,a,b})=0$.
In what follows we deal with the remaining case $q=p$.

From $v_{p}(k+1)\geq e_{p,k}$, we deduce that $p^{e_{p,k}}\mid
(k+1)$ and $p^e\mid (k+1)$ for each $e\in \{1,\ldots,e_{p,k}\}$. By
Lemma 3.6, any $p^e$ consecutive terms in the arithmetic progression
$\{b+ma\}_{m\in \mathbb{N}_0}$ are pairwise incongruent modulo $p^e$
since $p\nmid a$. Hence for each $e\in \{1,\ldots,e_{p,k}\}$, there
are exactly $\frac{k+1}{p^e}$ terms divisible by $p^e$ in any $k+1$
consecutive terms of the arithmetic progression $\{b+ma\}_{m\in
\mathbb{N}_0}$. So we have that $f_{e}(n)=\frac{k+1}{p^e}-1$ for
each $e\in \{1,\ldots,e_{p,k}\}$. In other words, for every $n\in
\mathbb{N}$, we have $f_{e}(n+1)=f_{e}(n)$. It then follows from
(3.2) that for every $n\in \mathbb{N}$, we have
$g_{p,k,a,b}(n+1)=g_{p,k,a,b}(n)$. Thus $P_{p,k,a,b}=1$ and
$v_{p}(P_{k,a,b})=0$. Therefore part (ii) is proved.
\end{proof}

\section{\bf Proof of Theorem 1.2}

In this section, we first prove Theorem 1.2.\\
\\
{\it Proof of Theorem 1.2:} By Theorem 2.5, we know that $g_{k,a,
b}$ is periodic. Denote by $P_{k, a, b}$ its smallest period. Let
first gcd$(a, b)=1$. Then by Lemma 3.7, for any prime $p$ such that
$p|a$, we have $v_p(P_{k,a,b})=0$. For any prime $p$ satisfying
$p\nmid a$ and $p\le k$, we have by Lemma 3.7,
$P_{p,k,a,b}=p^{v_p(P_{p,k,a,b})}=p^{v_p(P_{k,a,b})}$. So by
Proposition 3.8 we infer that
$$P_{k,a,b}=\prod_{p \ {\rm prime}, p\leq k}p^{e_{p}(k,a)},$$
where
\begin{align*}
e_{p}(k,a):={\left\{
  \begin{array}{rl}
0, \quad&\text{if} \ v_{p}(k+1)\geq e_{p, k} \ {\rm or} \ p\mid a,\\
e_{p, k}, \quad&\text{otherwise.}
 \end{array}
\right.}
\end{align*}
Using the integer $L_k$, we obtain immediately that $P_{k,a,b}=Q_{k,
a}$ as required, where $Q_{k, a}$ is defined as in (1.1).

Now let gcd$(a, b)>1$. If ${\rm gcd}(a,b)=d$ and $a=da'$ and
$b=db'$, then ${\rm gcd}(a',b')=1$ and we can easily check that
$g_{k,a,b}(n)=d^kg_{k,a',b'}(n)$ for any $n\in \mathbb{N}$. From
this one can easily derive that the periodic functions $g_{k,a,b}$
and $g_{k,a',b'}$ have the same smallest period, i.e., $
P_{k,a,b}=P_{k, a', b'}$.  But the result for the case gcd$(a, b)=1$
applied to the function $g_{k, a', b'}$ gives us that $P_{k, a',
b'}=Q_{k, a'}$, with $Q_{k, a'}$ defined as in (1.1). Thus the
desired result $P_{k,a,b}=Q_{k, a'}$ follows immediately. This
completes the proof of Theorem 1.2. \hfill$\Box$
\\

It was proved by Farhi and Kane \cite{[FK]} that there is at most
one prime $p\le k$ such that $v_{p}(k+1)\geq e_{p,k}$. We noticed
that such a prime $p$ was given in Proposition 3.3 of \cite{[FK]}
without the condition $p\le k$, but such restriction condition is
clearly necessary because otherwise Proposition 3.3 of \cite{[FK]}
would not be true. For example, letting $p$ be any prime number
greater than $k+1$ gives us $v_{p}(k+1)=0=e_{p,k}$. Comparing the
smallest period $P_{k,a,b}$ of the function $g_{k, a, b}$ with the
smallest period $P_k$ of the function $g_{k}=g_{k, 1, 0}$, we arrive
at the relation between $P_{k,a,b}$ and $P_k$ as follows:
$$
P_{k,a,b}=\frac{P_k}{\displaystyle\prod_{{\rm prime} \ p|{\rm
gcd}(a',P_k)}p^{e_{p,k}}},
$$
where $a'=\frac{a}{{\rm gcd}(a,b)}$. From this one can read that
$P_{k, a, b}=P_k$ if $a|b$.

Finally, we give an application of Theorem 1.2 as the conclusion of this paper.\\

\noindent{\bf Example 4.1.} Let us consider the least common
multiple of any $k+1$ consecutive positive odd numbers. To study
this problem, we define arithmetic function $h_{k}$ by
$$
h_{k}(n):=\frac{(2n+1)\cdot(2n+3)\cdots(2n+2k+1) } {{\rm
lcm}(2n+1,2n+3,\cdots,2n+2k+1)} \  \ (n\in \mathbb{N}).
$$
By Theorem 1.2, we know that $h_k$ is periodic and for any integer
$k\ge 2$, the exact period $R_k$ of $h_{k}$ is given by
$R_k=\frac{L_k}{2^{e_{2,k}}\cdot D_k}$, where
\begin{align*}
D_{k}={\left\{
  \begin{array}{rl}
p^{e_{p,k}},
 \quad& {\rm if} v_p(k+1)\geq e_{p, k} \
 {\rm for \ some \ odd \ prime} \ p\le k, \\
1, \quad&\text{otherwise.}
 \end{array}
\right.}
\end{align*}

\begin{center}
{\bf Acknowledgment}
\end{center}
The authors would like to thank Professor Smyth and the referees for
their careful reading of the manuscript and helpful suggestions
which improved its presentation.


\begin{thebibliography}{99}
\bibitem{[A]} T.M. Apostol,  {Introduction to analytic number
theory,} Springer-Verlag, New York, 1976.
\bibitem{[BK]} G. Bachman and T. Kessler, {On divisibility properties of certain multinomial
coefficients II,} J. Number Theory 106 (2004), 1-12.
\bibitem{[BKS]} P. Bateman, J. Kalb and A. Stenger, A limit involving
least common multiples, Amer. Math. Monthly 109 (2002), 393-394.
\bibitem {[BBGH]} M.A. Bennett, N. Bruin, K. Gy\"ory and L. Hajdu,
{Powers from products of consecutive terms in arithmetic
progression}, Proc. London Math. Soc. 92 (2006), 273-306.
\bibitem{[C]} J. Cilleruelo, The least common multiple of a quadratic sequence,
arXiv:1001.3438.
\bibitem{[F1]} B. Farhi, {Minoration non triviales du plus petit
commun multiple de certaines suites finies d'entiers,} C.R. Acad.
Sci. Paris, Ser. I 341 (2005), 469-474.
\bibitem{[F2]} B. Farhi, {Nontrivial lower bounds for the least common multiple of some
finite sequences of integers,} J. Number Theory 125 (2007), 393-411.
\bibitem{[FK]} B. Farhi and D. Kane, New results on the least common multiple of
consecutive integers, Proc. Amer. Math. Soc. 137 (2009), 1933-1939.
\bibitem{[GT]} B. Green and T. Tao, {The primes contain arbitrarily long arithmetic
progressions,} Ann. of Math. (2) 167 (2008), 481-547.
\bibitem{[Ha]} D. Hanson, {On the product of the primes,} Canad. Math.
Bull. 15 (1972), 33-37.
\bibitem{[HW]} G.H. Hardy and E.M. Wright, {An introduction
to the theory of numbers, Fourth Edition,} Oxford University Press,
London, 1960.
\bibitem{[HF]} S. Hong and W. Feng, {Lower bounds for the least common multiple
of finite arithmetic progressions,} C.R. Acad. Sci. Paris, Ser. I
343 (2006), 695-698.
\bibitem{[HK]} S. Hong and S.D. Kominers, Further improvements of lower bounds
for the least common multiple of arithmetic progressions, Proc.
Amer. Math. Soc. 138 (2010), 809-813.
\bibitem{[HLe]} S. Hong and K.S. Enoch Lee, Asymptotic behavior of
eigenvalues of reciprocal power LCM matrices, Glasgow Math. J. 50
(2008), 163-174.
\bibitem{[HLo]} S. Hong and R. Loewy, {Asymptotic behavior of eigenvalues
of greatest common divisor matrices,} Glasgow Math. J. 46 (2004),
551-569.
\bibitem{[HY1]} S. Hong and Y. Yang, {Improvements of lower bounds for the least
common multiple of finite arithmetic progressions,} Proc. Amer.
Math. Soc. 136 (2008), 4111-4114.
\bibitem{[HY2]} S. Hong and Y. Yang, On the periodicity of an arithmetical function, C.R. Acad.
Sci. Paris, Ser. I 346 (2008), 717-721.
\bibitem{[Hu]} L.-K. Hua, {Introduction to number theory,} Springer-Verlag, Berlin Heidelberg, 1982.
\bibitem{[MS]} G. Myerson and J. Sander, {What the least
common multiple divides II,} J. Number Theory 61 (1996), 67-84.
\bibitem{[N]} M. Nair, {On Chebyshev-type inequalities for primes,}
Amer. Math. Monthly 89 (1982), 126-129.
\bibitem {[S]} N. Saradha, {Squares in products with terms in
an arithmetic progression}, Acta Arith. 86 (1998), 27-43.
\bibitem {[SS]} N. Saradha and T.N. Shorey, {Almost squares in
arithmetic progression}, Compositio Math. 138 (2003), 73-111.

\end{thebibliography}
\end{document}